\documentclass[12pt]{article} 
\usepackage{amsmath}
\usepackage{amssymb} 
\usepackage{amscd}
\usepackage{amsthm}
\usepackage{enumitem}

\usepackage[utf8]{inputenc}
\usepackage{hyperref}

\newtheorem{theorem}{Theorem}

\newtheorem{lemma}[theorem]{Lemma}
\newtheorem{corollary}[theorem]{Corollary}

\def\dim{{\mbox{dim}}}

\def\Ext {{\mbox{Ext}}}

\def\Spin{{\mbox{Spin}}}

\def\calf{{\mathcal F}}

\def\bbbone{\mbox{\rm 1\hspace {-.6em} l}}

\begin{document}

\enlargethispage{3cm}

\thispagestyle{empty}
\begin{center}
{\bf CLIFFORD ALGEBRAS, MESON ALGEBRAS \\
\vspace{0.5cm}
AND HIGHER ORDER GENERALIZATIONS}\\
\vspace{0.2cm}

\end{center}
   
\vspace{1cm}

\begin{center}
Michel DUBOIS-VIOLETTE
\footnote{Laboratoire de Physique des 2 Infinis Ir\`ene Joliot Curie\\
	P\^ole Th\'eorie, 
	IJCLab UMR 9012\\
	CNRS, Universit\'e Paris-Saclay, 
	B\^atiment 210\\
	F-91406 Orsay Cedex, France\\
	michel.dubois-violette@ijclab.in2p3.fr} and 
Blas TORRECILLAS
\footnote {Universidad de Almeria\\ Departamento de Matem\'aticas\\ Carretera Sacramento s/n\\ 04120 La Ca\~nada de San Urbano, Almeria, Espa\~na\\
blas.torrecillas@gmail.com} \end{center}
\vspace{0,5cm}

\begin{abstract}
We analyse the homogeneous parts of Clifford and meson algebras and point out that for the Clifford algebra it is related to fermionic statistics, that is,  to fermionic parastatistics of order 1 while for the meson algebra it is related to fermionic parastatistics  of order 2. We extend these homogeneous algebras into corresponding algebras related to fermionic parastatistics of all orders. We then define correspondingly higher order generalizations of Clifford and meson algebras.
\end{abstract}

\vfill
\today

 \newpage
\tableofcontents

\newpage
\section{Introduction}

Throughout this paper $\mathbb K$ denotes a commutative field of characteristic zero and all vector spaces and algebras are over $\mathbb K$. The binomial coefficients are denoted by $C^p_N$, i.e. one has $(1+x)^N=\sum^N_{p=0} C^p_N x^p$. Concerning Young diagrams, Young tableaux and representation theory we refer to \cite{ful:1997}.\\

Given a finite-dimensional (pseudo-) euclidean vector space $E$ there are 3 closely related algebras \cite{jac:1968}: 

\begin{enumerate}[label=(\alph*)]
\item  The Jordan spin factor $JSpin(E)$\\
\item The Clifford algebra $C(E)$\\
\item The meson algebra $D(E)$.
\end{enumerate}
The {\sl Jordan spin factor} $JSpin(E)$ is the unital Jordan algebra $\mathbb K\bbbone \oplus E$ with unit $\bbbone$ and product $\circ$ defined by $x\circ y=(x,y)\bbbone$ for $x,y\in E$ where ($\cdot,\cdot $) denotes the scalar product.\\

The {\sl Clifford algebra} $C(E)$ is the unital associative algebra generated by the elements of $E$ with product satisfying 
\begin{equation}
xy+yx=2(x,y)\bbbone
\end{equation}
for any $x, y\in E$.\\
The {\sl meson algebra} $D(E)$ is the unital associative algebra generated by the elements of $E$ with product satisfying

\begin{equation}
xyx = (x,y) x
\end{equation}

for any $x, y \in E$.\\
In other words $C(E)$ and $D(E)$ are defined by

\begin{equation}
C(E)=T(E)/ (\{x\otimes y+y\otimes x - 2(x,y) \bbbone \vert \forall x,y\in E\})\tag{1'}
\end{equation}
and 
\begin{equation}
D(E)=T(E)/ (\{x\otimes y\otimes x-(x,y)x\vert \forall x,y\in E\})\tag{2'}
\end{equation}
where $T(E)$ is the tensor algebra of $E$ and where $(I)$ for $I \subset  T(E)$ denotes the 2-sided ideal of $T(E)$ generated by $I$.\\

$C(E)$ is the {\sl universal unital associative envelope} of the Jordan algebra $JSpin(E)$ while $D(E)$ is the {\sl universal unital multiplicative envelope} of the Jordan algebra $JSpin(E)$,\cite{jac:1968}. (This means that the left-modules over $D(E)$ are the Jordan modules over $JSpin(E))$ .\\

The fact that the canonical linear mapping of $JSpin(E)$ into $C(E)$ is injective means that the Jordan algebra $JSpin(E)$ is special.
The canonical mapping of $JSpin(E)$ into $D(E)$ is injective as the one of any Jordan algebra into its universal multiplicative envelope.\\

There is an injective homomorphism of unital algebras \linebreak[4] $\psi_2 : D(E)\rightarrow C(E)\otimes C(E)$ such that
\begin{equation}
\psi_2(x)=\frac{1}{2}(x\otimes \bbbone +\bbbone \otimes x)
\end{equation}
for  any $x\in E$,  and in fact $x\in J\Spin(E)$ $(JSpin(E)$ is  {\sl strongly special}), \cite{jac:1968}.\\

In this paper we shall be concerned with the {\sl homogeneous parts} of the above algebras and the generalizations of these algebras. The homogeneous parts of the Clifford algebra and of the meson algebra are also referred to as their {\sl neutral} versions since they are formally obtained by setting to zero the scalar product.\\

The structure of the neutral Clifford algebra $C_0(E)$ is without mystery since it coincides with the exterior algebra $\wedge E$ of $E$. The structure of the neutral meson algebra $D_0(E)$ is more subtle : It is the unital associative algebra generated by the elements of $E$ with product satisfying
\[
xyx=0
\]
for any $x,y\in E$. We shall study this algebra in detail and this will lead us to an interpretation in terms of a remarkable algebra connected with the fermionic parastatistics studied in \cite{mdv-pop:2002} and \cite{mdv-pop:2012c}. This is the key of our  higher order generalization.\\

The group of automorphisms of both the algebras $C(E)$ and $D(E)$ is the orthogonal group $O(E)$ of $E$ while the one of their homogeneous parts is the full linear group $GL(E)$ of $E$.\\

The algebras $C(E)$ and $D(E)$ are $\mathbb Z_2$-graded while their homogeneous parts are $\mathbb N$-graded and are isomorphic to their associated $\mathbb N$-graded algebras (PBW-property).\\
Moreover all these algebras are {\sl Frobenius algebras}.\\

The meson algebra is also referred to as the Duffin-Kemmer-Petiau algebra (DKP algebra). This algebra was introduced by G\'erard Petiau in 1936 and independently by R.J. Duffin in 1938 and by N. Kemmer in 1939, \cite{pet:1936}, \cite{duf:1938},\cite{kem:1939}.\\

For modern advances in meson algebras we refer to a recent article by J. Helmstetter \cite{hel:2024} and to its references.
 
\section{The neutral meson algebra $D_0(E)$}

As explained in the introduction the neutral meson algebra is the cubic algebra generated by the elements of $E$ with relations
\begin{equation}
xyx=0,\>\>\> \forall x, y \in E
\end{equation}
thus one has
\[
D_0(E)=T(E)/(\{x\otimes y \otimes x\vert x, y \in E\})
\]
where $T(E)$ is the tensor algebra of $E$.\\

Applying successively $x_1,\cdots, x_i,\cdots, x_n\in E$ , one has in  view of (4)  $x_i x_{i+1}x_{i+2}+x_{i+2} x_{i+1} x_i=0$. Thus the result in $D_0(E)$ is that one has complete antisymmetry between the $x_{2p}$ and between the $x_{2p+1}$ respectively. It follows that
\begin{equation}
D^n_0(E)=\left\{
\begin{array}{ll}
D^{2m}_0(E)=\wedge^m E\otimes \wedge^mE & \mbox{for}\>\> n=2m\\
D^{2m+1}_0(E)=\wedge^{m+1}E\otimes \wedge^m E &   \mbox{for}\>\> n=2m+1
\end{array}
\right.
\end{equation}
because (4) are the only constraints defining the associative unital algebra $D_0(E)$ which is $\mathbb N$-graded since it is homogeneous (i.e. cubic here). The Hilbert series $H(t)=\sum_n \dim(D^n_0(E))t^n$ of $D_0(E)$ reads :
\begin{equation}
H(t) =\sum_p(C^p_N)^2t^{2p}+ C^{p+1}_N C^p_N t^{2p+1}
\end{equation}
where $N=\dim(E)$. It follows that one has :
\begin{equation}
\dim(D_0(E))=\sum_p(C^p_N)^2+ C^p_N C^{p-1}_N=C^N_{2N+1}
\end{equation}
for the dimension of $D_0(E)$.\\

The counterpart of the injective homomorphism $\psi_2$ of (3) is the injective homomorphism $\varphi_2:D_0(E)\rightarrow \wedge E\otimes \wedge E$ such that
\begin{equation}\label{phi0}
\varphi_2(x)=\frac{1}{2}(x\otimes \bbbone + \bbbone \otimes x),\>\>\> \forall x\in E
\end{equation}
which is now invariant by the linear group $GL(E)$.\\

The decomposition of $D_0(E)$ into irreducible $GL(E)$-subspaces corresponds to non trivial combinatorial identities that we now describe.\\

The isomorphism classes of irreducible representations of $GL(E)$ are labeled by the Young diagrams $\lambda$ with at most $N$ lines (i.e. with columns of length $\leq N$).\\

Since for each degree $n$ $(0\leq n\leq 2N)$ $D^n_0(E)$ is invariant by $GL(E)$ one has to analyse $\wedge^m E \otimes \wedge^m E=D_0^{2m} (E)$ for $n$ even ($n=2m)$ and\linebreak[4] $\wedge^{m+1}E \otimes \wedge^mE=D^{2m+1}_0(E)$ for $n$ odd ($n=2m+1)$. In terms of Young diagram this corresponds to

\[
m\,  \left\{\,
 \begin{tabular}{|c|}\hline
\\  \hline \vdots
\\ \hline
\\ \hline
 \end{tabular} 
 \right.
  \otimes\,  m\,  \left\{\,
 \begin{tabular}{|c|}\hline
\\  \hline \vdots
\\ \hline
\\ \hline
 \end{tabular} 
 \right. \quad \text{for}\quad D^{2m}_0(E)\\
\]

 \text{and}\\
\[
 m+1\,  \left\{\,
 \begin{tabular}{|c|}\hline
\\  \hline \vdots 
\\ \hline
\\ \hline
\\\hline
 \end{tabular} 
 \right.
 \otimes\,  m\,  \left\{\,
 \begin{tabular}{|c|}\hline
\\  \hline \vdots
\\ \hline
\\ \hline
 \end{tabular} 
 \right.  \quad \text{for}\quad D^{2m+1}_0(E)
 \]
so one can obtain diagrams $\lambda$ with $2m$ boxes, $\vert\lambda\vert=2m$, and $2m+1$ boxes, $\vert\lambda\vert=2m+1$, having at most 2 columns \cite{ful:1997} of length $\leq N$. The dimensions  of the corresponding representations of $GL(E)$ is the number of tableaux with filling contained in $\{1,2,\cdots,N\}$ of the above diagrams. For such a 2-columns diagram $\lambda_{p,q}$ with columns of length $p,q$ with $N\geq p\geq q\geq 0$ and $p+q=\vert \lambda_{p,q}\vert$, this number, and therefore the dimension $\dim(\lambda_{p,q})$ of the representation of $GL(E)$, is given by
\begin{equation}
\dim(\lambda_{p,q})=\frac{p-q+1}{p+1} C^p_N C^q_{N+1}
\end{equation}
where the $C^p_N$ are the binomial coefficients.

\begin{lemma}\label{L1}
One has the following identity
\[
\frac{p-q+1}{p+1} C^p_N C^q_{N+1}=C^p_N C^q_N-C^{p+1}_N C^{q-1}_N (=\dim(\lambda_{p,q}))
\]
for $N\geq p\geq q\geq 0$.
\end{lemma}

\noindent\underbar{Proof}.   Indeed one has
\[
\begin{array}{lll}
\frac{p-q+1}{p+1} C^p_N C^q_{N+1}& =C^p_NC^q_{N+1}-\frac{q}{p+1} C^p_N C^q_{N+1}\\
\\
&=C^p_N C^q_{N+1}- C^{p+1}_{N+1}C^{q-1}_N\\ 
\\
&= C^p_N (C^q_N+C^{q-1}_N)-(C^{p+1}_N + C^p_N) C^{q-1}_N\\
\\
&= C^p_N C^q_N - C^{p+1}_N C^{q-1}_N. \square
\end{array}
\]

\begin{theorem} \label{T1}
The dimension $\dim(D^n_0(E))$ of the homogeneous part $D^n_0(E)$ of degree $n$ of the neutral meson algebra $D_0(E)$ satisfies the following identity

\[
\dim(D^n_0(E))= \sum_{\substack{{p+q=n}\\
{N\geq p\geq q\geq 0}}}  \dim (\lambda_{p,q})
\]
for any integer $n\leq 2N$.
\end{theorem}

\noindent\underbar{Proof}. For $n$ even say $n=2m$, one has in view of Lemma\ref{L1}

\[
\begin{array}{lll}
\sum_{k\geq 0} \dim (\lambda_{m+k,m-k}) &= (C^m_N)^2- C^{m+1}_N C^{m-1}_N\\
\\
& + (C^{m+1}_N  C^{m-1}_N - C^{m+2}_N C^{m-2}_N)+ \cdots\\
\\
& + (C^{m+k}_N C^{m-k}_N - C^{m+k+1}_N C^{m-k-1})+ \cdots
\end{array}
\]
Thus finally $\sum_{k\geq 0} \dim (\lambda_{m+k,m-k})= (C^m_N)^2=\dim(D^{2m}_0(E))$ since all other terms cancel $(C^{m+k}_N C^{m-k}_N$\, \text{for}\,  $k\geq 1$).\\

For $n$ odd, $n=2m+1$\\
\[
\begin{array}{lll}
\sum_{k\geq 0}\dim(\lambda_{m+1+k,m-k}) &=C^{m+1}_N C^m_N - C^{m+2}_N C^{m-1}_N \\
\\
&+ (C^{m+2}_N C^{m-1}_N - C^{m+3}_N C^{m-2}_N) + \cdots
\end{array}
\]

$\Rightarrow \sum_{k\geq 0} \dim (\lambda_{m+1+k,m-k})=C^{m+1}_N C^m_N=D^{2m+1}_0(E). \square$\\

Last theorem means that the representation of $GL(E)$ in $D^n_0(E)$ is the direct sum of the irreductible representations corresponding to the Young diagrams with two columns $\lambda_{p,q}$ with $N\geq p\geq q\geq 0$ and $p+q=n$, (each occurring with multiplicity one).  The degree in the 3-homogeneous algebra $D_0(E)$ being obviously invariant by the action of $GL(E)$ on $D_0(E)$, this is equivalent to the fact that the representation of $GL(E)$ in $D_0(E)$ is the direct sum of the irreducible representations corresponding to the 2-columns Young diagrams i.e. one has
\begin{equation}
\dim(D_0(E))=\sum_{N\geq p\geq q\geq 0}\dim(\lambda_{p,q})
\end{equation}
in terms of dimensions.

\section{Another approach and its generalization}
\subsection{The key result}
The key of our generalization is the following theorem.
\begin{theorem}\label{T2}
One has the equivalence
\[
xyx=0\quad \forall x,y\in E\Longleftrightarrow \left\{\begin{array}{l}
[[x,y],z]=0\quad \forall x,y,z\in E\\
\mbox{and}\\
x^3=0\quad \forall x\in E
\end{array}
\right.
\]
in any associative algebra generated by $E$.
\end{theorem}

\noindent\underbar{Proof}. $xyx=0 \quad \forall x,y\in E$ implies clearly
\[
[[x,y],z]=(xyz+ zyx) - (zxy + yxz)=0
\]
and $x^3=0 \quad \forall x,y,z \in E$.\\
Conversely, the identity
\[
3 xyx=[[x,y],x] + xyx + x^2y + yx^2
\]
implies that $xyx=0$ if $[[x,y],z]=0$ and $x^3=0$ $\forall x,y,z\in E$ since $xyx+x^2y +y x^2$ is the linear part in $y$ of $(x+y)^3$. $\square$

\begin{corollary}\label{C1}
The homogeneous meson algebra $D_0(E)$ is the quotient of 
\[
\calf(E)=T(E)/(\{[[x,y],z] \vert x,y,z\in E\})
\] 
by the two sided ideal generated by the $x^3$ for $x\in E$ i.e. one has :
\[
D_0(E)=\calf(E)/(\{x^3\vert x\in E\}).
\]
\end{corollary}

\subsection{The remarkable algebra $\calf(E)$}
The relations $[[x,y],z]=0$ for $x,y,z\in E$ of $\calf(E)$ are the relations of the creation (resp. destruction) operators of a general parafermi system with $N=\dim(E)$ degrees of freedom. Whenever one adds the relations $x^{n+1}=0$ to these relations one obtains the relations of the creation (resp. destruction) operators of a parafermi system
of order $n$ (order of the parafermi system),\cite{gre:1953}, \cite{gre:1966}, \cite{ohn-kam:1982} . 

Thus from the above results $D_0(E)$ describes the algebra of creation (resp. destruction) operators of a parafermi system of order 2 while $C_0(E)=\wedge E$ corresponds to a parafermi system of order 1 that is of an ordinary Fermi system. Thus Formula (8) corresponds  then to the ``Green ansatz" (or Green representation) for $n=2$. Let us now describe the two main properties of $\calf(E)$.\\

The first property is that $\calf(E)$ is the universal enveloping algebra of the free 2-nilpotent Lie algebra generated by $E$, that is the graded Lie algebra $E\oplus\wedge^2E$ with Lie bracket $[x,y]=x\wedge y$ if $x,y\in E$ and $[x,y]=0$ otherwise. It follows that the graded algebra $\calf(E)$ is an Artin-Schelter algebra of global dimension $N(N+1)/2=\dim(E+\wedge^2E)$\cite{flo-vat:2011}
and that its Yoneda algebra $\Ext_\calf(\mathbb K, \mathbb K)$ is the cohomology $ H(E\oplus \wedge^2 E)$ the Lie algebra $E\oplus \wedge^2 E$, that is one has
\[
\Ext_\calf(\mathbb K, \mathbb K)=\oplus_i \Ext^i_\calf (\mathbb K,\mathbb K)
\]
with $\Ext^i_\calf(\mathbb K,\mathbb K)=H^i(E\oplus \wedge^2 E)$ and in fact the Yoneda algebra of $\calf(E)$ is bigraded
\[
\Ext_\calf(\mathbb K, \mathbb K)=\oplus_{i,j} \Ext^{i,j}_\calf (\mathbb K, \mathbb K)
\]
where $i$ is the homological degree while $j$ is the degree induced by the tensorial degree in $E$.\\

The second property of $\calf(E)$ is its content in representations of $GL(E)$. Indeed,  by construction, the linear group $GL(E)$ acts by polynomial representations on $\calf(E)$ and it turns out that each irreducible representation appears here with multiplicity one. Thus by identifying an irreducible representation of $GL(E)$ with the corresponding Young diagram, one can write for $\calf(E)$ the following expansion
\begin{equation}
\calf(E)=\oplus \lambda
\end{equation}
where the summation runs over all Young diagrams $\lambda$ with columns of length (number of boxes) smaller or equal to $N=\dim (E)$, while, for its Yoneda algebra one can write
\begin{equation}
\Ext _\calf(\mathbb K, \mathbb K)= \oplus_{\lambda=\bar\lambda} \lambda
\end{equation}
where $\lambda$ runs over all symmetric diagrams with columns of length smaller or equal to $N$, ($\lambda\mapsto \bar\lambda$ denotes the symmetry with respect to the main diagonal).

\subsection{Higher order generalizations}
The above analysis shows that there is an obvious generalization of the neutral Clifford algebra $C_0(E)=\wedge E$ and the neutral meson algebra $D_0(E)$ for higher order $n$ greater than 2. Namely let $\Phi_n(E)$ be the quotient of $\calf(E)$ by the two-sided ideal generated by the subset $\{x^{n+1}\vert x\in E\}$ of $\calf(E)$, that is the unital algebra
\begin{equation}
\Phi_n(E)= \calf(E)/(\{x^{n+1}\vert x\in E\})
\end{equation}
for any integer $n\geq 1$. Then $\Phi_1(E)=\wedge E$ is the neutral Clifford algebra $C_0(E)$ and $\Phi_2(E)=D_0(E)$ is the neutral meson algebra while, more generally, $\Phi_n(E)$ corresponds to the fermionic parastatistics of order $n$ (see 3.2). Furthermore the content of $\Phi_n(E)$ in irreducible representations of $GL(E)$ consists in term of Young diagrams of all Young diagrams with at most n columns and $N=\dim(E)$ lines, each occurring with multiplicity one.\\

\noindent One has an injective algebra homomorphism $\varphi_n:\Phi_n(E)\longrightarrow  \underbrace{\wedge E \otimes \cdots \otimes \wedge E}_{n}$ induced by  
\begin{equation}\label{eq14}
\varphi_n(x) = \frac{1}{n} (x\otimes \bbbone^{\otimes^{n-1}} + \bbbone \otimes x\otimes \bbbone^{\otimes^{n-2}}+ \cdots + \bbbone^{\otimes^{n-1}}\otimes x)
\end{equation}
for $x\in E$ which generalizes (\ref{phi0}) for arbitrary $n\geq 1$. The homomorphism $\varphi_n$ describes the ``Green ansatz" \cite{gre:1953},\cite{gre:1966}, \cite{ohn-kam:1982} .
\\
\\
For $n=3$ for instance, $x^4=0$ $\forall x\in E$ and (\ref{eq14}) reads
\[
\varphi_3(x) =\frac{1}{3}(x\otimes \bbbone \otimes \bbbone + \bbbone \otimes x\otimes \bbbone+\bbbone \otimes \bbbone \otimes x)
\]
for the homomorphism $\varphi_3$.  
\section{Higher order generalizations of Clifford and meson algebras}

Now the way to define higher order generalizations of Clifford and meson algebras is straightforward. We first notice that the Green ansatz for order $n$ for the neutral version is the injective homomorphism $\varphi_n:\Phi_n(E)\rightarrow \underbrace{\wedge E\otimes \cdots \otimes \wedge E}_n$ given by (\ref{eq14}). Then we define $\Psi_n(E)$ to be the subalgebra of $\underbrace{C(E)\otimes\cdots\otimes C(E)}_n$, generated by  the
\begin{equation}\label{eq15}
\psi_n(x)=\frac{1}{n}(x\otimes \bbbone^{\otimes^{n-1}}+\cdots +\bbbone^{\otimes^{n-1}}\otimes x)
\end{equation}
for $x\in E$.
It is clear that $\Psi_1(E)=C(E)$ and that $\Psi_2(E)=D(E)$ are the Clifford and the meson algebras while (\ref{eq15}) is the version of (\ref{eq14}), for $(\cdot,\cdot) \not = 0$.\\

Let us remind that the relations which define the neutral version $\Phi_n(E)$ are $[[x,y],z]=0$ and $x^{n+1}=0$ for $x,y,z\in E$. Let us study what happens to these terms $[[x,y],z]$ and $x^{n+1}$ in $\Psi_n(E)$.\\

For $[[x,y],z]$ the result is given by the following lemma.
\begin{lemma}\label{L2}
In $\Psi_n(E)$ one has
\[
[[x,y],z]=\frac{4}{n^2}((y,z)x-(x,z)y)
\]
for $x,y,z\in E\subset \Psi_n(E)$.
\end{lemma}

\noindent \underbar{Proof}.
For $[[x,y],z]$ one has in the Clifford algebra $C(E)=\Psi_1(E)$
\begin{equation}\label{eq16}
[[x,y],z]=4((y,z)x-(x,z)y)
\end{equation}
by using the definition of $C(E)$. It then follows easily that in $\Psi_n(E)$ one has
\begin{equation}\label{eq17}
[[x,y],z]=\frac{4}{n^2}((y,z)x-(x,z)y)
\end{equation}
$\forall x,y, z\in E$. Indeed one verifies that 
\[
[[\psi_n(x), \psi_n(y)], \psi_n(z)]=\frac{4}{n^2} ((y,z)\psi_n(x)-(x,z) \psi_n(y))
\]
by using the definition of $\psi_n(x)$ for $x\in E$ in Formula (\ref{eq15}).$\square$\\

For $x^{n+1}$ the computation is more delicate, in any case, since the right hand side vanishes for the neutral version, one has a development
\begin{equation}\label{18}
x^{n+1}=\sum^{p\leq \frac{n+1}{2}}_{p=1} a_p\parallel x \parallel^{2p} x^{n-2p+1}
\end{equation}
with coefficients $a_p\in \mathbb K$. The computation is not a priori a problem but  it is not obvious how to give a general form for the coefficients $a_p$.\\

Of course for $n=1$ and $n=2$ i.e. for the Clifford and the meson algebras there is only one coefficient which is $1\in \mathbb K$, that is $a_1=1$ ; for the Clifford algebra $x^2=\parallel x\parallel^2\bbbone$ while for the meson algebra $x^3=\parallel x\parallel^2 x$.\\

Let us show the results of the computations up to $n=7$.\\

For $n=1$ and $n=2$ one has
\[
x^2=\parallel x\parallel^2\>\> \mbox{or}\>\> x^2-\parallel x \parallel^2=0
\]
\[
x^3=\parallel x\parallel^2 x\>\> \mbox{or}\>\>( x^2-\parallel x\parallel^2)x =0
\]
for $n=3$, 4, 5, 6,7 one finds
\[
x^4=\frac{10}{3^2}\parallel x\parallel^2x^2-\frac{1}{3^2}\parallel x\parallel^4
\]
\[
x^5=\frac{5}{4}\parallel x\parallel^2x^3-\frac{1}{4}\parallel x\parallel^4 x
\]
\[
x^6=\frac{7}{5}\parallel x\parallel^2x^4-\frac{259}{5^4}\parallel x\parallel^4x^2+\frac{3^2}{5^4}\parallel x\parallel^6
\]
\[
x^7=\frac{14}{3^2}\parallel x\parallel^2 x^5-\frac{7^2}{3^4}\parallel x\parallel^4x^3+\frac{2}{3^4}\parallel x\parallel^6 x
\]
\[
x^8=\frac{12}{7}\parallel x\parallel^2 x^6-\frac{282}{7^3}\parallel x\parallel^4x^4+\frac{12916}{7^6}\parallel x\parallel^6 x^2-\frac{3^2 5^2}{7^6}\parallel x\parallel^8
\]
Moreover all these results factorize as
\[
(x^2-(\frac{1}{3})^2\parallel x\parallel^2)(x^2-\parallel x\parallel^2)=0
\]
\[
x(x^2-(\frac{1}{2})^2\parallel x\parallel^2)(x^2-\parallel x\parallel^2)=0
\]
\[
(x^2-(\frac{1}{5})^2\parallel x\parallel^2)(x^2-(\frac{3}{5})^2\parallel x\parallel^2) (x^2-\parallel x\parallel^2)=0
\]
\[
x(x^2-(\frac{1}{3})^2\parallel x\parallel^2)(x^2-(\frac{2}{3})^2\parallel x\parallel^2)(x^2-\parallel x\parallel^2)=0
\]
\[
(x^2-(\frac{1}{7})^2\parallel x\parallel^2)(x^2-(\frac{3}{7})^2 \parallel x\parallel^2)(x^2-(\frac{5}{7})^2\parallel x\parallel^2))(x^2-\parallel x\parallel^2)=0
\]
Therefore this suggests the following for $n=2p+1$ and for $n=2p+2$, $(p\in \mathbb N)$

\begin{equation}\label{eq19}
\prod^p_{q=0}\left(x^2-\left(\frac{2q+1}{2p+1}\right)^2\parallel x\parallel^2\right)=0
\end{equation}
 for $n=2p+1$ and
 \begin{equation}\label{eq20}
 x\prod^p_{q=0}\left(x^2-\left(\frac{2q+2}{2p+2}\right)^2\parallel x\parallel^2\right) =0
 \end{equation}
 
that is
\begin{equation}
x\prod^p_{q=0} \left( x^2- \left(\frac{q+1}{p+1}\right)^2\parallel x\parallel^2\right) =0\tag{20'}
\end{equation}
for $n=2p+2$, ($p\in \mathbb N$).\\

Let us introduce $X_n=n\psi_n(x)\in C(E)^{\otimes^n}$ that is
\begin{equation}\label{eq21}
X_n=x\otimes \bbbone^{\otimes^{n-1}}+\bbbone \otimes x\otimes \bbbone^{\otimes^{n-2}}+\cdots + \bbbone^{\otimes^{n-1}}\otimes x
\end{equation}
for $x\in E\subset C(E)$ then (\ref{eq19}) and (\ref{eq20}) read
\begin{equation}\label{eq22}
\prod^p_{q=0}\left(X^2_{2p+1}-(2q+1)^2\parallel x\parallel^2\right)=0
\end{equation}
for $n=2p+1$ and
\begin{equation}\label{eq23}
X_{2p+2}\prod^p_{q=0}\left( X^2_{2p+2}-(2q+2)^2\parallel x\parallel^2\right)=0
\end{equation}
for $n=2p+2$.
\begin{lemma}\label{L3}
Relations (\ref{eq22}) are satisfied in $C(E)^{\otimes^{2p+1}}$ and relations (\ref{eq23}) are satisfied in $C(E)^{\otimes^{2p+2}}$ for $p\in \mathbb N$.
\end{lemma}

\noindent \underbar{Proof}. 
Let us start with $n$ odd i.e. $n=2p+1$ ($p\in \mathbb N$). For $x\in E\subset C(E)$, one has $x^2=\parallel x\parallel^2$. On the other hand the $2p+1$ terms of the left-hand side of (\ref{eq21}) $x\otimes \bbbone^{2p}$, $\bbbone\otimes x\otimes \bbbone^{2p-1},\cdots, \bbbone^{\otimes 2p}\otimes x$ commutes so can be simultaneously diagonalised.\\

Assume first that the $x's$ are all in the same eigen value say $x=\parallel x\parallel$. Then the corresponding eigenvalue of $X_{2p+1}$ is $(2p+1)\parallel x\parallel$ and of course if $x=-\parallel x\parallel$ (the other eigenvalue)one has $X_{2p+1}=-(2p+1)\parallel x\parallel$.  So in this case one has $X^2_{2p+1}-(2p+1)^2\parallel x\parallel ^2=0$.\\

Assume now that all the $x's$ have the same eigenvalue $\parallel x\parallel$ for instance except one which has the opposite value $-\parallel x\parallel$. It follows that $X_{2p+1}=(2p+1)\parallel x\parallel - 2 \parallel x\parallel =(2p-1)\parallel x\parallel$.Therefore one has $X^2_{2p+1}-(2p-1)^2\parallel x\parallel^2=0$ by the same argument as above.\\

Thus by assuming now that $p-q$ of the $x's$ has an eigenvalue opposite than the other ones for $p\geq q$ then one gets $X^2_{2p+1}-(2q+1)^2\parallel x\parallel^2=0$. \\

Thus by taking all possibilities one gets
\[
\prod^p_{q=0}(X^2_{2p+1}- (2q+1)^2 \parallel x\parallel^2)=0
\]
that is Relation (\ref{eq22}).\\
Similarily, by the same argument one gets Relation (\ref{eq23}) for the even case $n=2p+2$. $\square$

The following theorem is a corollary of Lemma (\ref{L2}) and Lemma (\ref{L3}) above.
\begin{theorem} \label{FD}
The algebra $\Psi_{2p+1}(E)$ is the unital associative algebra generated by $E$ with relations (odd case)
\[
\left\{
\begin{array}{l}
\displaystyle{[[x,y],z]=\left(\frac{2}{2p+1}\right)^2((y,z)x-(x,z)y)\>\>\>  \forall x,y,z\in E}\\
\\
\displaystyle{\prod^{q=p}_{q=0}}\left(x^2-\left (\frac{2q+1}{2p+1}\right)^2\parallel x\parallel^2\right)=0 \>\>\> \forall x\in E
\end{array}
\right.
\]
while the algebra $\Psi_{2p+2}(E)$ is the unital associative algebra generated by $E$ with relations (even case)
\[
\left\{
\begin{array}{l}
\displaystyle{[[x,y],z]=\left(\frac{1}{p+1}\right)^2((y,z)x-(x,z)y) \>\>\> \forall x,y,z\in E}\\
\\
\displaystyle{x\prod^{q=p}_{q=0} \left( x^2 - \left(\frac{q+1}{p+1}\right)^2 \parallel x\parallel^2\right)=0\>\>\> \forall x\in E}
\end{array}
\right.
\]
for any $p\in \mathbb N$.
\end{theorem}

In fact $\Psi_n(E)$ is a deformation of its neutral version $\Phi_n(E)$.

\section{Conclusion}
We have described generalization in higher order of the Clifford and the meson algebras and of their neutral versions. As pointed out in Section 3 these homogeneous versions are directly related to fermionic parastatistics of order 1 and order 2 respectively. Thus in our generalization the order means fermionic parastatistics order so it is worth noticing here that this order has nothing to do with the order considered in the interesting paper \cite{hel-mic-rev:2012}. In fact the connection between the Clifford algebra and the fermionic statistics is well known but, as pointed out here, the meson algebra has the same connection with the fermionic parastatistics of order 2, (fermionic statistic = fermionic parastatistics of order 1).\\
Now as pointed in the introduction both the Clifford algebra $C(E)=\Psi_1(E)$ and the meson algebra $D(E)=\Psi_2(E)$ are the solution of universal problems concerning the Jordan algebra $JSpin(E) $, \cite{jac:1968} so it is natural to investigate whether there is something similar connecting each $\Psi_n(E)$ with the Jordan spin factor $JSpin(E)$ for $n>2$. More generally it remains to study in details the properties of the algebras $\Psi_n(E)$ and of their neutral versions $\Phi_n(E)$.\\

\noindent {\bf Acknowledgements}

\medskip

Work supported by the project Lanzadera and PID2024-158993NB-100.
Part of this work was done while the second author was visiting IHES. The author wishes to thank the institution for its hospitality and support.
\newpage 


\end{document}